\documentclass[10pt,a4paper]{article}
\usepackage{t1enc}
\usepackage[latin1]{inputenc}
\usepackage[english]{babel}
\usepackage{harvard}
\begin{document}
\thispagestyle{empty}
\begin{center}
\subsection*{Philosophy as a Cultural Resource and Medium of Reflection for Hermann Weyl}
{\em Erhard Scholz, Wuppertal}\footnote{I want to thank the organizers of the conference at Cret Berard for their interest and the opportunity to present the remarks on H. Weyl to a predominantly philosophical audience,  John Stachel and Sk\'uli Sigurdsson for discussions, linguistic support and critical remarks to earlier versions of the paper.}

\end{center}
\begin{abstract}
 Here we review a kind of  post-World-War-II "Nachtrag" to H. Weyl's philosophical comments on mathematics and the natural sciences published in the middle of the 1920s. In a talk given at Z\"urich in the late 1940s, Weyl discussed F.Gonseth's dialectical epistemology and considered it  as being restricted too strictly to aspects of  historical change. His own experiences with
 post-Kantian dialectical philosophy, in particular J.G. Fichte's derivation of the concept of space and matter, had been a stronger dialectical  background for his own 1918 studies in purely infintitesimal geometry and the early geometrically unified field theory of matter (extending the Mie-Hilbert program). Although now Weyl distantiated himself  from the speculative features of his youthful  philosophizing and in particular from his earlier enthusiasm for Fichte, he again had deep doubts as to the cultural foundations of modern mathematical sciences and its role in material culture of high modernity.  For Weyl,  philosophical "reflection" was a cultural necessity; he now turned towards K. Jasper's and M. Heidegger's existentialism to find deeper grounds, similar to his turn towards Fichte's philosophy after World War I.

\end{abstract}

\subsubsection*{Introduction}
Philosophers  might get  interested in  Hermann Weyl for   many  reasons.    One of these may  be  his explicitly  philosophical publications, first of all and best known  his {\em Philosophie der Mathematik und Naturwissenschaften}
\cite{Weyl:PMN} written in 1926 as a contribution to the {\em Handbuch der Philosophie} (English  translation with comments and appendices 1949). Moreover his mathematical and scientific work  contains in many places, sometimes in decisive branching regions of his thought, philosophical arguments from which Weyl drew motivations for the scientific questions  posed or the research orientation chosen.

Among philosophies to which Weyl referred, Husserl's phenomenology is the best known. Weyl's relationship to phenomenology  has been  investigated and documented in several publications,\footnote{Among them \cite{vanDalen:Husserl}, \cite{Tonietti:Deppert}, \cite{daSilva:Weyl}, \cite{Tieszen:Weyl} and \cite{Mancosu/Ryckman}. Da Silva argues convincingly for a close methodological link between Husserl and Weyl's {\em Das Kontinuum}, whereas he is much more cautious on  Husserl's overall role for Weyl than Tonietti.  I tend to be even more cautious and  to agree with  the sceptical remarks in this respect by  Paolo Mancosu and Thomas Ryckman in \cite{Mancosu/Ryckman} at the occasion of their  edition of and comments on the correspondence between Weyl and Oskar Becker. See below (beginning of the section  ``Fichte's dialectical construction\ldots'').}
but Weyl never became a devoted adherent of any single philosophy.
He rather was a wanderer through the philosophical fields differing with the changes in his scientific views,\footnote{An early example his alienation from a (narrow) Kantianism under the influence of his awareness of Hilbert's work on the foundations of geometry at the beginning of his university studies in 1904/05} his personal environment,\footnote{His personal relationship with Husserl's  student Helene Joseph, beginning in  1912 and soon to become his wife, and his acquaintance  and later friendship with Fritz Medicus, professor  of pedagogy and philosophy in Z\"urich, after 1913 had strong influences on the philosophical world of discourse, in which Weyl chose his own point of view.} and with  cultural and social breaks in his environment. In particular he strongly felt the influences of both great wars of the 
20th century and was deeply affected by the social and cultural changes resulting from them.

On this occasion I do not  want to go into either Weyl's more systematic exposition of  a {\em Philosophy of Mathematics and Natural Sciences} \cite{Weyl:PMN2} or his relationship with Husserl's philosophy.
Rather I shall discuss some aspects of  his  relationship with dialectical philosophies  of his time, those of  Ferdinand Gonseth and Martin Heidegger,   or of some identifiable influence on his scientific work, like
Johann Gottlieb Fichte's {\em Wissenschaftslehre}, which Weyl came close to in  his early years at Z\"urich, under the  influence of  Fritz  Medicus.\footnote{For better evaluation of the character of the ``influence'' see the discussion in the last section of this paper.}  

\subsubsection*{Some  remarks by Weyl on Gonseth's dialectical philosophy and on Heidegger}
In his  philosophically and foundationally (with respect to mathematics) most radical phase, between   1916 and 1923/24, Weyl  pursued  his intention of basing mathematics and physics on ``immediate intuition (unmittelbare Einsicht)'' \cite[146]{Weyl:Krise} or on what he called an ``analysis of essence attempted by phenomenological philosophy (von der ph\"anomenologischen Philosophie (Hus\-serl) angestrebte Wesensanalyse)'' \cite[$^4$1921, 133]{Weyl:RZM}. At that time he rejected the more formal methods of construction of mathematical knowledge proposed, among others, by his former teacher David Hilbert. In the middle of the 1920s he reapproached Hilbert's position on the foundations of mathematics to a certain degree and conceded to formal mathematics  founded on axiomatic principles  a potentiality for   the intellectual appropriation  of external reality or its symbolic representation, as Weyl preferred to formulate it,  which went far beyond what he had accepted some years earlier. In  a phase of restructuring of his thought between  1920 and 1922, he came to the conviction that some ``transcendent reality'' (God, matter \dots) lies beyond the realm of phenomenologically accessible reality,   which at least in some aspects   (e.g. matter understood as a dynamical agent) might be indirectly represented by mathematical symbols the nature and exact form of which had still to be explored (for a short phase in the early 1920s he  considered singularities of classical fields as a possible candidate) or  invented and developed (at the end of the 1920s he turned towards spinor fields with an additional $U(1)$ gauge symmetry and, at least in principle, their mathematically highly problematical second quantization). In this case he realized  that formal mathematics might even have an advantage over   immediately insightful (``phenomenological'') mathematics, because  in its conceptual constitution it was free  from the restrictions of the latter. I will call this view the {\em symbolic realism} of the ``mature''  Weyl.

The topic of symbolic realism was taken up and elaborated by Weyl in his contribution   ``Wissenschaft als symbolische Konstruktion des Menschen (Science as Symbolic Construction of Man)''  to the {\em Eranos Jahrbuch}  \cite{Weyl:142}, in which he sketched, now as a sexagenarian, a kind of r\'esum\'e  of the transition from the early modern,  mechanistic image of nature and the role of mathematics in it to the state it had assumed in the classical modern view of the first half of the 20th century. This r\'esum\'e  is comparable to and continues the one  he had attempted in his {\em Philosophie der Mathematik und der Naturwissenschaften} more than twenty years earlier. It was quite different in content from the first one, due to  further advances, not  at all reassuring for Weyl,  in the foundations of mathematics and in fundamental physics. 

Weyl characterized the quintessence of the transition as a replacement  of the old ``mechanistic construction of the world'' in the sense of some  assumed spatio-temporal materiality understood  in an atomistic sense by a ``construction in pure symbols'' \cite[295]{Weyl:142}. He deplored the fact that physicists and philosophers had both continued to \ldots
\begin{quote}
\ldots stick stubbornly to the principles of a mechanistic interpretation of the world (mechanistische Weltdeutung) after physics had, in its factual structure, already outgrown the latter. They have the same excuse as the inhabitant of the mainland (Landmensch) who for the first time travels on the open sea: he will desparately try to stay in sight of the vanishing coast line, as long as there is no other coast in sight, towards which he steers. \cite[299]{Weyl:142}
\end{quote}
In his own perspective of a ``construction of the world in pure symbols'',  his symbolic realism as I called  it above, Weyl hoped to see the first contours of a ``new coast line''. He frankly admitted, however, that he was not able to exclude the possibility that, perhaps,  
 ``\ldots we are only deceived by a fog-bank'' (ibid.).
On that background he discussed Einstein's {\em principle of relativity} (approvingly) and  Bohr's {\em principle of complementarity}, which he wanted to accept only as far as it ``corresponds to a fact in quantum physics, which can be given a precise mathematical meaning'' \cite[338]{Weyl:142}. We should take into account for this passage that in 1948 about twenty years of hard work had been invested by physicists and mathematicians to link quantum physics with special relativity. About the time of Weyl's talk the young generation of theoretical physicists were just achieving a decisive step forward with renormalizable quantum electrodynamics \cite{Schweber:QED}, whereas the wider perspective of a synthesis of gravitation and quantum physics was still deeply covered by some barely identifiable ``fog-bank''.

At this place, Weyl made a short excursion  to the outlook of Gonseth's philosophy, as presented, e.g. in \cite{Gonseth1943}:
\begin{quote}
The call for a dialectization of knowledge rings out from Z\"urich. It is not completely clear to me what this means. \cite[339]{Weyl:142}
 \end{quote}
He  agreed with  Gonseth that  modern science continues to create a theoretical image of the world, although  this image and the language in which it is formulated  is no longer restricted by any fixed a priori but is shaped by scientific experience and  `` its interpretations (Deutungen) which are forced upon it by the context of  science itself''. Although Gonseth had, of course, analyzed the epistemic role of  some of these {\em Deutungen} in a stronger sense than Weyl presented in this short reference, i.e. as  ``constitutive schemata''  in knowledge  production,\footnote{Cf. \cite[here, pp. 89ff.]{Heinzmann:Gonseth}.}  Weyl essentially agreed with Gonseth thus far.  But, in the evaluation of the range of the change from the formerly fixed a  priori to the constructed constitutive schemata,  he had a difference that he frankly stated: \begin{quote}
But this free outlook, which does justice to the interaction between construction and reflection, probably does  not yet deserve the label of dialectic. (ibid.)
\end{quote}
Such a critical remark was not at all dismissive with respect to dialectical thought, as one might   expect from a mathematical scientist active in  one of  the  ``two cultures'' of the 
20th century.\footnote{Cf. \cite{Skuli:DMV}.}
 As distinct from most of the contemporary scientists and mathematicians, Weyl had become acquainted with the  dialectical philosophy of (post-Kantian) German idealism in his early Z\"urich years under the influence of F. Medicus, and had,  in particular,  intensely studied Fichte's {\em Wissenschafts\-lehre}. As a result, he requested from  a proper dialectic that it should be richer in the interplay of internal oppositions than just  allowing for conceptual shifts in the  historical development of knowledge systems. 
Following a short presentation of Gonseth's discussion of the concept of  velocity and its modification by relativity, he commented:
\begin{quote}
In the new relativistic picture the original concepts are `lifted (aufgehoben)' in  Hegel's double meaning of the word. That may be accurate in the historical sense, but even then it would only be a `historical' dialectic. (ibid.) 
\end{quote}
In the late 1940s, Weyl had distanced himself considerably   from Fichte. As he explicitly stated, he now sympathized with another version of dialectical thought,   ``existential philosophy, e.g. in the Heideggerian form'' \cite[343]{Weyl:142}. To hint, at least at the direction in which he looked for a connection between existentialist philosophy and modern science, Weyl presented a simple but strong consideration that  would counteract or  even prevent any attempt to strive for a new reductionistic world picture.

He reminded his audience that   quantum physics (QP) posed a much more fundamental riddle than relativity theory (RT). In RT the classical concepts of kinematics and dynamics were ``lifted'' in a way that allowed a reconstruction and understanding of prerelativistic thought well defined mathematically and physically. Critically  reflecting  the Copenhaguen interpretation of quantum mechanics from his point of view,  he contrasted this well defined conceptual connection in RT with the dilemma outlined by Bohr for  QP, which arose from the necessity to represent  the ``hidden physical process'', by a mathematical symbolism dissociated from classical physics and everyday experience. On the other hand  he declared that  ``natural understanding of the world and the language in which it is expressed'' is necessary for the process of measurement and observation and  its description, ``perhaps slightly purified and clarified (gereinigt und gekl\"art) by classical physics''. He saw no possibility for a complete theory of measurement and observation in terms of QP. On the assumption that such a difference might continue to hold  in the future, he  commented:
\begin{quote}
Then we had, here,   a true dialectic, impossible to lift by any historical development: the hidden  soil (der dunkle Boden)  [of everyday life and natural language, ES] can be illuminated from the higher viewpoint  of quantum physics, but it remains the ground (Grund) that cannot be dissolved into the light of the higher region. \cite[341]{Weyl:142} 
\end{quote}
 Such an  unresolvable opposition reminded Weyl of a comparable one in existential philosophy. There the {\em Dasein} (being-there) as the form of being, which includes
 self-awareness, is achieved by human beings rising to   consciousness of themselves.  On the other hand the same process necessarily led to a ``burial of the external world in nothingness'', and thus made it necessary later to ``prove'' the latter indirectly. The rise of human cognition to a rationalist
self-consciousness  of an  individual (Cartesian) self,  led necessarily to  what Weyl now called the ``problem of the external world (Problem der Aussenwelt)''.\footnote{This formulation is close in language and content to the ``problem of matter (Problem der Materie)'' that Weyl had encountered in the later 1920s, after he had given up the reductionist idea of deriving matter structures, in a Mie-type approach, from pure (classical) field theory.} 
\begin{quote}
(O)ne attempts to glue the isolated subject left behind to the disconnected patches of the world, but it remains  a patchwork. \cite[344]{Weyl:142}]
\end{quote} 
Although Weyl admitted that his argument for comparability between the different oppositions in QP and existential  philosophy  had no  compelling  power, he nevertheless insisted on the value of the problem posed by this observation. Moreover, he indicated how the existential split between the clarity of self-conscious being-there (Dasein) and the dark patches of insight of the external world might be resolved, not by a progression in pure thought, but by a transition to  another way of life and being, in which  the self might be percieved as part of a communicating web of 
self-conscious beings and the broader world in which they live.  
\begin{quote}
In understanding myself as being-together-with (seiend-mit)  I also understand different
 being-there (Dasein). This is not, however, a knowledge achieved and concluded by new cognition, but a primary existential way of being that is the {\em conditio sine qua non} for cognition and  knowledge. (ibid.) 
\end{quote}
We find other notes and speeches which show that, after the experience of the destructive powers unleashed by the Second World War,  existential philosophy  had become   for Weyl a medium of  reflection not necessarily leading to consequences for scientific knowledge, but endowed with a   value of its own. Thirty years earlier in his life, in his encounter with Fichte's philosophy during and shortly after the First World War,  things had been slightly different in this respect. Weyl had  been attracted by philosophical considerations as a young man,  but at that time  he  tried to connect his newly-gained philosophical convictions to his ongoing scientific research.

 \subsubsection*{Fichte's dialectical construction of the concept of space}
In his glance back at his intellectual life, {\em Erkenntnis und Besinnung} \cite{Weyl:166}, Weyl recounted his turn towards Fichte's {\em Wissenschaftslehre} under the influence of F. Medicus. He admitted to having been seized by Fichte's idealist metaphysics, although he  had been lucky enough to find a counterweight in Helene's (his wife's) more sober style of thought, shaped by  Husserl's phenomenology: 
\begin{quote}
I  had to concede to her that Fichte, by his stubbornness in pursuing an idea, blind to nature and facts,  was swept away into increasingly abstruse constructions. \cite[637]{Weyl:166}  
\end{quote}
On the other hand, he had been impressed and attracted by Fichte, whom in 1954 he still esteemed  as a
\begin{quote}
\ldots constructivist of the purest water, who pursues his independent path of conctruction without looking right or left \cite[641]{Weyl:166}.
\end{quote}

This characterization gives a clue as to why the young  mathematician Weyl was attracted by Fichte in the years between 1916 and 1921.  At that time Weyl did not refer much  to Fichte  in  his publications; nevertheless he did not completely hide the attraction and sympathy he felt for Fichte's philosophy. On the first pages of {\em Das Kontinuum}, he discussed the difficult  problem  of a philosophically satisfying clarification of the foundations of logic and confessed:
\begin{quote}
We cannot hope to give here a final clarification of the essence of fact, judgement, object, property; this task leads into metaphysical abysses; about these one has to seek advice from men whose name  cannot be stated without  earning a compassionate smile (mitleidiges L\"acheln) --- e.g. Fichte.   \cite[2]{Weyl:Kont}
\end{quote}
At the time  more esteem could  be expected among mathematicians, in particular from the G\"ottingen milieu, from references to Husserls's philosophy.  Thus we find more frequently public references to Husserl's {\em Ph\"anomenologie} than to Fichte's {\em Wissenschaftslehre}, initially restricted, however, to the  discussion of the problem of time and Husserl's {\em Ideen zu einer reinen Ph\"anomenologie} \cite{Husserl:Ideen}.\footnote{References  in {\em Das Kontinuum} \cite[iv]{Weyl:Kont} and at  the beginning of  {\em Raum - Zeit - Materie} (Space - Time - Matter) \cite[4, footn. 1]{Weyl:RZM}.}   {\em After} Weyl had formulated his ``purely infinitesimal'' geometry, he looked for a more philosophical underpinning of his new geometry, and proceeded, in 1920 and the following years,  to an a priori derivation of its structure by what he then called  the mathematical analysis of the problem of space ({\em mathematische Analyse des Raumproblems}). His turn of thought coincided with the publication of Husserl's altered  (second) edition of {\em Logische Untersuchungen} in the years 1921 (vol. II) and 1922 (vol. I). P. Mancosu and T. Ryckman argue that  Weyl {\em only then} started to read the latter with some care. In fact,  in  a
 letter-of-thanks  (March 26/27, 1921) for the gift of the 2nd edition of the {\em Logische Untersuchungen} by its author, Weyl wrote to Husserl that he had ``only now'' been able to get some ``superficial'' acquaintance  with the book's  content \cite[290]{Husserl:Korr_7}.     

When he  included a first sketch of his investigations on the {\em Raumproblem} in the 4th edition of {\em Raum - Zeit - Materie} it appeared fitting to him to declare his recent investigations to be an 
``analysis of essence''  as attempted by phenomenological philosophy (``von der ph\"anomenologischen Philosophie (Husserl) angestrebte Wesensanalyse'')'' \cite[$^4$1921, 133]{Weyl:RZM}. This remark included, now, a  global and non-explicit reference to Husserl's {\em Logische Untersuchungen} \cite{Husserl:LU}. We better abstain, however, from reading into it a long-range impact of Husserl's phenomenology on Weyl's own researches.\footnote{The dialogue between Weyl and Husserl may also have found  expression the other way round. In the 2nd edition of vol. I of {\em Logische Untersuchungen}, we find interesting discussions of the {\em Idee der reinen Mannigfaltigkeitslehre} and the division of labour among mathematicians/scientists and philosophers \cite[vol. I ($^2$1921), 258ff.]{Husserl:LU}.}
 His letter to Husserl was written three to four years {\em after} Weyl's first fundamental contributions to the foundations  of mathematics and ``purely infinitesimal geometry'' (gauge geometry) in 1918, nearly two years after his ``conversion'' to Brouwer's intuitionism, and at a time when he was coming close to the end of his most radical phase in foundations of mathematics and  unified field theory based on a dynamistic natural philosophy. At this time  the ``mathematical analysis of the problem of space''   opened up for him  a path towards a continuation of his  interest in gauge geometrical structures without  a necessary  link to a purely field theoretic (Mie type)  explanation of matter.
Even taking into account the appropriate modesty of Weyl as a scientific correspondent with respect to the professional philosopher, we find here a convincing additional evidence for Weyl's own later evaluation that, although Husserl's phenomenological philosophy had helped  him in finding a path towards   a ``freer view of the world''  than the positivist one, it had not been  the strongest  intellectual reference for him during  his years of {\em Sturm and Drang}  \cite[637]{Weyl:166}. 

As I have already argued on another occasion \cite{Scholz:ICM}, Fichte's dialectical construction of the concept of space seems to have been 
 of more importance for the development of Weyl's research orientation in this early, philosophically radical, phase of his work. In order to make the argument more accessible (and thus criticizable) I want  to give a short outline of some features of Fichte's philosophy as a ``radical constructivist''  at work in (dialectically) constructing the concept of space, both topics close to Weyl's mind in the years between  1916 and 1920.\footnote{See also the remarks on Fichte in (da Silva 1997).}

Fichte considered dialectics as the form of movement in the self-construction of thought, and described this conception in words that came to meet Weyl's needs in his search for intellectual orientation and support in the ``metaphysical abysses''  of the philosophical clarification of the foundations of logic. Fichte explained, e.g., in this context that he was not keen to develop the ``dialectics of making up or inventing something (die Dialektik des Aus- oder Erdenkens)''; he wanted to achieve a way of organizing the thought process in a way such that ``we are seized by  evidence''. This was, according to Fichte, the task of a ``true dialectic''.
\begin{quote}
Ingenuity gives only sudden evidence which can vanish again; true dialectic is the lawful method to achieve evidence. \ldots {\em Construction} is  now   the instruction to invent the concept by the imaginative power,  such that evidence may be gained. \cite[188]{Fichte:EinlWiss}\footnote{``Dies ist nicht die Dialektik des Aus- und Erdenkens, sondern das Denken macht sich uns selbst, die Evidenz ergreift uns. Durch Genie nur pl\"otzliche Evidenz, die wieder entschwinden kann; wahre Dialektik aber die gesetzm\"a\ss ige Methode, zu dieser Evidenz zu kommen. Die Kunst der Dialektik, wie alle Kunst, ist unendlich; nicht aber die Wahrheit. {\em Construction} ist nun die Anleitung, durch die Einbildungkraft den Begriff zu erfinden, da\ss{} die Evidenz sich einstelle. Es giebt da gleichsam eine urspr\"ungliche, und nach dieser ringen wir. Wie in der Mathematik: nicht die Wahrheit wird gemacht, sondern nur der Vortrag der Wahrheit.'' (emphasis in the original)}
\end{quote}  
A comparable claim for a direct link between construction and evidence  runs like a continuous thread through Weyl's radical writings on the foundations of analysis; in fact, on the level of  basic methodological convictions, it  unites the different, more technical positions which he adhered to between 1916 and 1920/21, from the  predicative, arithmetico-logical foundation in \cite{Weyl:Kont} to the
 (quasi-)intuitionistic strategy of \cite{Weyl:Krise}.

As a philosophical resource of inspiration for Weyl,  Fichte's construction or ``derivation'' of the concept of space seems to have been of far larger range. I sketch the most important  features (in the Weylian context) following Fichte's exposition in the {\em Grundri\ss{}   des Eigenth\"umlichen der Wissenschaftslehre} \cite{Fichte:Grundriss}. 

Fichte started from two opposing, initially separated, ``external'' perceptions, i.e. perceptions of appearances, $X$ and $Y$ (Fichte's own symbolism,  \cite[194]{Fichte:Grundriss}). As appearances  $X$ and $Y$ were, for Fichte, expressions of two  {\em forces} that had to be assumed ``necessarily'' (ibid. 197). Space was then developed by Fichte as the ``initially still unknown'' external determination of the possible relations between $X$ and $Y$  (i.e., in particular not generated by the interior determinations of the appearances $X$ and $Y$). For this determination (which later in the development --- or ``construction'' 
---  would become space) he introduced the formal symbol $O$. 

Now, the ``mode how $X$ is thus determined'', i.e. by its relation to the totalitiy of possible relations, may be symbolized as $x$, and the ``mode how $Y$ is thus determined'' as $y$. In this way, Fichte  introduced  symbols of external positioning $x$ and $y$ for the appearances  $X$ and $Y$,\footnote{In fact, Fichte used different small letters,   $z$ and $v$, in place of our  $x$ and $y$, which we use  here for mnemonic reasons.} 
and emphasized that the relation between $X$ and $Y$, and thus $O$ (the concept of space in generation), could ``by no means be derived from the self, but must be attributed to the things themselves'' (ibid. 195). 

Fichte  continued his argumentat that $X$ and $Y$ suffer ``exclusion and continuity (Ausschlie\ss ung und Continuit\"at)'' and that by such an interplay a ``common sphere'' is generated.\footnote{``Es wird demnach durch absolute Spontaneit\"at der Einbildungskraft eine solche gemeinschaftliche Sph\"are producirt'' (ibid. 196).}

 For the concept of space in generation Fichte postulated:
\begin{quote}
$O$ must be something which leaves the liberty of both ($X$ and $Y$, ES) in their activity completely undisturbed. \cite[198]{Fichte:Grundriss}\footnote{To readers acquainted with Weyl's analysis of the problem of space, a  striking similiarity to the latter's ``postulate of freedom'' may come to mind, when reading this Fichtean step.}
\end{quote}
Correspondingly, the external positionings $x$ and $y$ then represent  also ``spheres of activity (Sph\"aren der Wirksamkeit)'' of the forces ``assumed necessarily''  for the appearances  $X$ and $Y$. These must not, initially,  be considered from an extensional point of view, i.e. spacelike. Only after the formation of a ``common sphere'' of both
\begin{quote}
 [$O$ will be ] \ldots  posited (gesezt) as extended, connected, infinitely divisible and is {\em Space}. (ibid. 200).
\end{quote}
 Once the  concept of space has been derived, the positions ($x$ and $y$ etc.) are to be understood as result of an infinite division of space, as `ìnfinitely smallest parts of space''. Fichte was quite definite, though, that such parts could not be considered as points:
\begin{quote}
The infinitely smallest part of space is always a space, something endowed with continuity, not at all a mere point or the boundary between specified places in space. (ibid., 200)
\end{quote}

In these infinitesimal parts\footnote{``Infinitesimal'' is a word introduced by me (ES) not by Fichte, although the concept is.} 
the ``imaginative power (Einbildungskraft)''\footnote{Probably of the ``absolute self''.}
puts up a force which expresses itself ``with necessity''. Fichte thus arrived at a common conceptualization of space and systems of forces that went far beyond the Kantian dynamistic concept of matter in the {\em Metaphysische Anfangsgr\"unde der Naturwissenschaften}.
\begin{quote}
Thus intensity (here: force, ES) and extensity (space, ES) are by necessity synthetically united, and one must not try to discuss the one without the other. \cite[201]{Fichte:Grundriss}
\end{quote}
Fichte concluded:
\begin{quote}
Each force, by its necessary product, fills (\ldots) a place in space; and space is nothing but that which is filled or to be filled by these products. (ibid, 201) 
\end{quote}
I have  already said that ``places in space (Stellen im Raume)'' were not to be considered as points but as ``infinitely smallest parts''.

To summarize, Fichte ``constructed'' a concept of space starting from ``spheres of activities'' ($x$, $y$, etc.) of appearance ($X$, $Y$, \ldots); the system $O$ of possible external relations of these was ``posited'' as ``continuous, connected, and infinitely divisible''. Integration  of the conceptual structure transformed the spheres of activities ($x$, $y$, \ldots) into infinitesimal parts of $O$, which thus became the symbol and concept of space. The infinitesimal parts were  carriers of forces which Fichte considered as the primordial space-filling structures.

In  slightly more mathematical language we find here a clear and beautiful affinity to three
essential topics in Weyl's ``purely infinitesimal geometry'' and his first unified gauge
 field-and-matter theory of 1918:\footnote{On this background it does not come as a great surprise that T. Tonietti, in his excursion into Weyl's ``purely infinitesimal'' geometry, realizes at the height of his exposition of this point that he apparently  lost contact with Husserl, and that  he had better come back to the question of {\em Zeitbewu\ss{}tsein}, in order to  `` \ldots again encounter Husserl who otherwise might be forgotten'' \cite[366]{Tonietti:Deppert}.}
\begin{itemize}
\item[---] construction of a spacelike continuum from infinitesimal parts (not as a set of  points, endowed only in an additional step with a continuous, differentiable, etc. structure, as intended --- not yet achieved ---  by mainstream mathematics of the time),
\item[---] characterization of the space-filling entities as forces, the actions of which were initially specified only in the infinitesimal parts (Weyl's postulate to build geometrical structures ``purely infinitesimally''), 
\item[---] formation of matter as a form of appearance of space-filling forces (dynamical theory of matter, going back to Kant and mathematically rejuvenated by Mie and Hilbert).
\end{itemize}

\subsubsection*{Weyl's ``purely infinitesimal geometry'' and  his field/matter theory of 1918}
In order not to be misunderstood, I want to state the obvious in advance: Weyl's purely infinitesimal geometry and (gauge) field theory of 1918 were, of course, much  more than  a mere transformation of Fichtean ideas into mathematics and physics. There was a  substantial and complex background in the mathematical and physical knowledge of the time, 
deeply shaped by the work of Albert Einstein and  extended by 
G. Mie, D. Hilbert, T. Levi-Civita, and others.\footnote{For a comprehensive view, see \cite{Vizgin:UFT}, \cite{Stachel:HR}, \cite{Cao:FT}, for Mie and Hilbert \cite{Corry:Hilbert_Red,Corry:Hilbert_Mie}. }
 The intellectual work to be done was too demanding and technically too cumbersome, to allow such a (bad) reduction. In addition, our mathematician was confronted with the beneficial warnings of Helene Weyl not to follow   Fichte's  ``obstinacy, blind to nature and facts,'' in pursuing his ideas,  without sufficient reflection and intellectual distance. The over-enthusiastic language in key publications of the years 1918 to 1921 shows that he was, in fact, in such a danger. He experienced these years of social and political disorder as a time of  cultural crisis, and   was all the more longing for a  path towards an alleged transcendent order.\footnote{Cf. \cite{Skuli:DMV}.}

Already  by November 1917, when he finished  his work on  {\em Das Kontinuum} with its predicative construction of the arithmetically definable subset of the classical real continuum, Weyl was dissatisfied with his own theory.\footnote{For Weyl's 1918 approach to foundations of analysis see \cite{Feferman:Weyl2000} and  \cite[315ff.]{Cole/Korte},  for a discussion of his shifts in foundations \cite{Scholz:Hendricks}.} In public he compared it with Husserl's discussion of the continuum of time, which had a ``non-atomistic'' character, in drastic contrast to his own theory. No determined point of time can be exhibited, only  approximative fixing is possible, just as is the case for ``the continuum of spatial intuition (Kontinuum der r\"aumlichen Anschauung)'' \cite[70f.]{Weyl:Kont}.  Fichte had insisted on the same feature, as Weyl knew well at the time, but preferred not to touch upon explicitly at this place. However, he himself also accepted  the necessity that  the mathematical concept of continuum, the {\em continuous manifold}, should {\em not} be characterized in terms of modern set theory enriched by topological axioms (and axioms for coordinate systems), because this would contradict the concept (the ``essence'' as Weyl liked to claim) of the continuum. He continued to insist on this methodological principle even after softening his radical position or even turning away  from it. In a paper written in 1925, although published only posthumously, we  find, perhaps, the clearest expression of this Weylian conviction with respect to the set-theoretic approach to manifolds:
\begin{quote}
It seems  clear that it (set theory, ES) violates against  the essence of the continuum, which, by its very nature, can not at all be battered into a set of single elements. Not the relationship of an element to a set, but of a part to the whole ought to be taken as a basis for the analysis of the continuum. \cite[5]{Weyl:RGI}
\end{quote}
For Weyl, single points of a continuum were empty abstractions; an acceptable abstraction would arise only from a limiting procedure of localization inside a continuum:
\begin{quote}
The mf. is continuous, if the points are fused together in such a manner that it is impossible to display a single point, but always only together with a vaguely limited surrounding halo (Hof), with a `neighbourhood'. (ibid., 2)
\end{quote}
With such a conception Weyl entered difficult terrain, as no mathematical conceptual frame was in sight, which could satisfy  his  methodological postulate in a sufficiently elaborate way. For some years Weyl sympathized with Brouwer's idea to characterize points in the intuitionistic one-dimensional continuum by ``free choice sequences'' of nested  intervals  \cite{Weyl:Krise}.\footnote{More details of this aspect are given in \cite{Scholz:Hendricks}} 
He even tried to extend the idea to higher dimensions and explored  the possibility of a purely combinatorial approach to the concept of manifold, in which point-like localizations were given only by infinite sequences of nested star neighbourhoods in barycentric subdivisions  of a combinatorially defined ``manifold''. There arose, however, the problem of how to characterize the local manifold property in purely combinatorial terms. Although Weyl outlined a strategy to, possibly, overcome  this problem, his approach remained essentially without long-range effects in 20th century topology.

Weyl was much more successful on another level of investigation, in his attempts to rebuild differential geometry  in manifolds from a  ``purely  infinitesimal''  point of view. He generalized Riemann's proposal for a differential geometric metric 
\[ ds^2(x) = \sum_{i,j =1}^{n} g_{ij}(x)dx^i dx^j\]   
From his  purely infinitesimal point of view,   it seemed a strange effect that the length of two vectors $\xi (x)$ and $\eta (x')$ given at different points $x$ and $x'$ can be immediately (and in this sense ``objectively'') be compared in this framework after the calculation of
\[   |\xi (x)|^2 =  \sum_{i,j =1}^{n} g_{ij}(x) \xi ^i \xi ^j \, , \; \;    |\eta  (x')|^2 =  \sum_{i,j =1}^{n} g_{ij}(x') \eta  ^i \eta  ^j  \]

In this  context, it was comparatively easy, for Weyl, to give a purely infinitesimal characterization  of metrical concepts. He started from the well known structure of a conformal metric, i.e. an equivalence class $[g]$ of semi-Riemannian metrics $g = (g_{ij}(x))$ and $g' = (g_{ij}'(x))$ which are equal up to a point dependent positive factor $\lambda (x) > 0$, $g' = \lambda g$. Then, comparison of length made immediate sense only for vectors attached to the same ``point'' $x$, independently of the ``gauge'' of the metric, i.e. the choice of representative in the conformal class. To achieve comparability of lengths of vectors  {\em inside each infinitesimal neighbourhood} Weyl invented the concept of  {\em length connection} formed in analogy to the affine connection $\Gamma $, just distilled from the classical Christoffel symbols $\Gamma _{ij}^k $ of Riemannian geometry by Levi Civita (and himself, as far as the abstraction from the metrical context is concerned).\footnote{Cf. \cite{Reich:Levi}, \cite{Botta:Levi}.}
 The localization inside such an infinitesimal neighbourhood was given, as  would have been done already by the mathematicians of the past, by coordinate parameters  $x$ and $x' = x + dx$ for some ``infinitesimal  displacement'' $dx$.  Weyl's  length connection consisted, then, in an equivalence class of differential 1-forms $[\varphi]$, $\varphi = \sum_{i=1}^{n} \varphi _i dx^i$, where an equivalent representation of the form is given by 
$\varphi ' = \varphi -d \log \lambda  $, corresponding to a change of gauge of the conformal metric by the factor $\lambda $. Weyl called this transformation, which he recognized as necessary for the consistency of his extended symbol system, the {\em gauge transformation} of the length connection. That was the first step towards the investigation of (generalized) connections in modern differential geometry, i.e. no longer necessarily in the linear group and derived from a Riemannian metric. 

Weyl  thus established     a ``purely infinitesimal'' gauge geometry, where lengths of vectors (or derived metrical concepts in tensor fields) were immediately comparable only in the infinitesimal neighbourhood of one point, and for points of ``finite'' distance only after an integration procedure. This integration turned out to be, in general,   path dependent. Independence of the choice of path between two points $x$ and $x'$ holds if and only if  the ``length curvature'' vanishes. The concept of curvature was built in direct analogy to the curvature of the affine connection and turned out to be, in this case, just the exterior derivative of the length connection: $f = d \varphi$. This led Weyl to a coherent and conceptually pleasing realization of a metrical differential geometry built upon ``purely infinitesimal'' principles. He explored some features of this geometry, which appeared of immediate importance to him: existence and uniqueness of a compatible affine connection, gauge behaviour (gauge invariance respectively  ``gauge covariance'' specified by ``gauge weight'')  of tensors derived from the data $[g], \; [\varphi]$ of the Weylian metric, relationship to Riemannian manifolds, in particular different natural choices of gauge in the (semi-)Riemannian case, etc.

Moreover, Weyl was convinced of important consequences of his new gauge geometry for physics. The infininitesimal neighbourhoods understood as ``spheres of activity'', as Fichte might have said, suggested looking for interpretations of the length connection as a ``field''  representing physically active quantities. And in fact, building on the (mathematically obvious) observation $df = 0$, which  (in coordinates) was formally identical with the second system of the generally covariant Maxwell equations, Weyl immediately drew the conclusion that the length curvature $f$ ought to be identified with the electromagnetic field. For Weyl this had some, initially equally, important consequences: 
\begin{itemize}
\item[(1)] Conservation of current appeared as a consequence of gauge invariance.
\item[(2)] Conformal metric, the gauge principle, length connection and affine connection were part of a unified conceptually convincing structure.
\item[(3)] Gravitation (with its potential $g$) and electromagnetism (potential $\varphi$) were intrinsically unified and both identified with constitutive elements of the geometric structure.
\item[(4)] The Mie-Hilbert theory of a combined Lagrange function $L(g,  \varphi)$ for the action of  the gravitational field ($g$) and  electromagnetism  ($\varphi$) was further geometrized and technically enriched by the  principle of gauge invariance (for $L$).
\end{itemize}
Weyl thus believed for a short while (from 1918 to  summer 1920) that his gauge geometric modification of the Mie-Hilbert theory would finally show the path towards  the blue flower of dynamism, i.e. an  explanation of the stable structures of the elementary particles known at the time (electron and proton), the mathematical derivation of the basic stable composite configuration (atoms), at least in the simplest cases, and their energy levels (spectra). 

He gave up the belief, however, in the ontological correctness of the purely field-theoretic approach to matter (point (4) above)  in late summer 1920,  substituting in its place a philosophically motivated a priori argumentation for the conceptual superiority of his gauge geometry (i.e. point (2)) during the following years (his analysis of the problem of space, 1921 --- 1923).\footnote{Cf.  \cite[215ff.]{Cole/Korte} and \cite[85ff.]{Scholz:RZM}.} At the end of the 1920s  he finally withdrew the metrical version of the gauge characterization of electromagnetism (point (3) above) and transformed it into the new form of a $U(1)$ gauge connection in the  extended symmetries of a  generally convariant formulation of the Dirac equation, which he proposed in 1929.\footnote{Cf. \cite{Straumann:DMV} and \cite[chap.5]{Vizgin:UFT}.} The gauge  interpretation of current conservation (point (1)), finally, was upheld in the new form. In the early 1930s Pauli assimilated it  to the knowledge of quantum physics, and it became, in a generalized form, a central methodological  element in the self-construction of theoretical fundamental physics during the second half of the 20th century.\footnote{Cf. \cite{ORaif/Strau}.} 

The goal of a unified description of gravitation and electromagnetism, and the derivation of matter structures from it, was nothing specific to Weyl. After 1915 it had become a hot topic  in G\"ottingen, after Hilbert's incursion into general relativity.\footnote{Cf. \cite{Vizgin:UFT}, \cite{Corry:Hilbert_Mie}, \cite{Corry:Hilbert_Red}, \cite{Renn/Stachel:Hilbert}, \cite{Rowe:HilbertGRT}, \cite{Goldstein/Ritter:UFT}.}
Hilbert had built upon Mie's work, and Mie was fond of Kant's dynamistic approach to the problem of matter and other purely field-theoretic theories of the electron mass from the turn of the century. So Weyl's attempt at a unification was  in its general perspective not at all foreign to the larger enterprise of contemporary mathematical physics and natural philosophy, at least in the scientific environment of G\"ottingen. In Weyl's theory, the ``purely infinitesimal'' approach to manifolds and the ensuing possibility to geometrically unify the two known interaction fields gravitation and electromagnetism, this perspective took on  a particularly dense and conceptually sophisticated shape.  Thus it may be not surprising, although perhaps a bit ironic, 
 that Hilbert showed awareness of the ``extreme idealization'' of the geometric classical field theoretic approach to fundamental physics only after he could reflect on Weyl's sharpening of their common approach in the winter 1919/20 \cite[98ff.]{Hilbert:Rowe}. 
Hilbert's surprise regarding  Weyl's theoretical program which was so close  in its goals  to his own enterprise, only conceptually more densely knit,  helped him, apparently,  to reflect critically  on the  attempts to derive matter structures from classical fields. He now felt, that the ``completed process of physical idealization'' building upon classically deterministic mathematical structures would necessarily lead to what he now called  pejoratively ``Hegelian physics'' \cite[100]{Hilbert:Rowe}, as a label  for an all-embracing complete determinism. Such a world view would necessarily  run into the paradoxes of  time inversion also  holding for organic processes and, even stranger, the following effect:
\begin{quote}
Decisions in the proper sense could not take place at all, and the whole world process would not transcend the {\em limited content of a finite thought.} \ldots For the same reason also {\em intellectual culture} (das {\em Geistige}), in particular our thought, had to be something merely ephemeral (etwas blo\ss{} Scheinbares) --- an absurd consequence for a view of nature which results from the attempt to make all elements of reality (all Inhalte der Wirklichkeit) accessible to our thought. \cite[100f., emph. in original]{Hilbert:Rowe} 
\end{quote}

Hilbert was thus struggling, like Weyl, with the problem of how classical determinism in the theory of nature could go in hand with a unifying world view. In his stubborn rejection of  
post-Kantian dialectical philosophy, expressed in his misled stereotype of ``Hegelian physics'',  he pointed to a field that had served as a cultural background for Weyl. He apparently was not aware that  in its concrete  manifestation it was  Johann Gottlieb Fichte's  philosophy and not Hegel's, which had inspired the scientific imagination of his former student.\footnote{If Hilbert had known, he probably would not have cared.}
Moreover, his  polemics contain an unnoticed and  unvoluntary self-ironic twist:
Hilbert apparently did not realize that the verdict of a ``Hegelian physics'' might just as well, or even with stronger justification than to Weyl's program, be applied to his own theory.\footnote{J\"urgen Renn and John Stachel thus call Hilbert's original program  of 1915ff. an attempt of a ``theory of everything'' \cite{Renn/Stachel:Hilbert}.}

\subsubsection*{Weyl's late reflections on Heidegger}
What Weyl owed to Fichte's philosophy  cannot and should not be  described as   a ``philosophical  influence'' on mathematics, not even in the slightly more restricted version of a presumed ``influence'' of Fichte's philosophy of space and force on differential geometry, foundations of mathematics, and/or field physics. Fichte's philosophy was, of course, no historical agent or even actor. Its role in the genesis of Weyl's contribution to infinitesimal geometry, the concept of continuum, and field theory is much more precisely described as that of a {\em cultural resource}  which Weyl utilized. This perspective appears helpful for leaving behind the disciplinary boundaries in the investigation of   the history of mathematics and science in general.\footnote{I thank  Sk\'uli Sigurdsson who in our discussions insisted on such a distinction and emphasized, in particular, the   ``resource aspect'' of  philosophy  for Weyl; compare \cite{Skuli:DMV}}
Not in all cases, however, could Weyl's resort to  philosophers or philosophies be adequately described as a resource made serviceable in his scientific work. We have seen above, how Weyl referred to Heidegger's existentialist or Gonseth's dialectical philosophies in 1948. He was fond of, and probably needed, a broad reflection on his activities  as part of a conscious embedding in a cultural setting. Weyl considered such reflections as one part of a dyad on the other side of wich he saw and practised scientific activity, ``construction'' as he liked to say. Although he  talked about them in terms of  ``Besinnung'' (contemplation), which sounds surprisingly inoffensive  and
harmonious, these reflections  were not at all of  a harmonizing and introspective  nature. 

Just  to the contrary. Weyl was deeply bewildered by the great crises of the 20th century, the First World War and the following revolutionary disturbances, the rise  of Nazism in Germany, the Second World War with the holocaust and  the ``(let us hope short-lived) Nazification of the European continent'',\footnote{H. Weyl to  A. Johnson, March 22, 1941, quoted from  
\cite[285]{SiegmundSch:Rock}.}
and the ensuing menace of nuclear destruction. The goal of his ``Besinnung'' was to cope with such disconcerting experiences and to find some defendable position towards them. Weyl's late turn towards  Heidegger's philosophy is apparently an expression of the  experience  of being ``cast'' into a world  that had turned out to be socio-culturally unreliable and had developed strong features of destructiveness. Perhaps the strongest expression of his deep cultural doubts at this time can be found in a manuscript on {\em The development of mathematics since 1900}, in which he drew a parallel between Aristotle's characterization of metaphysics and of his own thoughts about modern mathematics:
\begin{quote}
\ldots Here some words of Aristotle come to my mind which, to be true,
refer to metaphysics rather than mathematics. Stressing its uselessness
as much as Hardy does in his apology of  mathematics, but at
the same time  its divinity, he says (Metaphysics 982b):
`For this reason its acquisition might justly be supposed to be beyond
human power; since in many aspects human nature is servile;
in which case, as Simonides says, `God alone can have this privilege',
and man should only seek the knowledge which is of concern to him
($\tau \grave{\eta} \nu$  $\kappa  \alpha \vartheta $'  
$\alpha \dot{\upsilon} \tau \grave{o}\nu$ $\dot{\epsilon} \pi \iota \sigma \tau \acute{\eta} \mu \eta \nu $). 
Indeed if the prets are right and the Deity is by nature jealous,
it is probable that in this case they would be particularly jealous and all
those who step beyond
($\pi \acute{\alpha} \nu \tau \alpha \varsigma \;  \tau o \grave{\upsilon} \varsigma  \; 
\pi \epsilon \rho \iota  \tau \tau o \acute{\upsilon} \varsigma  $)
are liable to misfortune.'

I am not so sure whether we mathematicians during the last decades have
not  `stepped beyond' the human realm by our abstractions. Aristotle, who
actually speaks about metaphysics rather than mathematics, comforts us by
hinting that the envy of the Gods is but a lie of the prets
(`prets tell many a lie' as the proverb says).

For us today the idea that the Gods from which we wrestled the
secret of knowledge by symbolic construction will revenge our
$ \upsilon \beta \rho \iota \varsigma $ has taken on a quite concrete form.
For who can close his eyes against the menace of our self-destruction
by science; the alarming fact is that the rapid progress of scientific
knowledge is unparalleled by a congruous growth of man's
moral strength and responsibility, which has hardly chance in
historical time. \cite[6f.]{Weyl:Math20Jhdt} 
\end{quote} 
Thus Weyl's ``Besinnung'' did not at all  avoid or compensate for uncomfortable experiences. He posed and presented  problems arising in mathematics, culture and society as sharply as he could. Although he feared the menace  arising from what he, in his educated language shaped by the classical humanist tradition,   called the ``hubris'' of modern mankind, he did not give up the hope that another mode of (existential) being might be possible. In his formulation inspired by Heidegger, cited at the beginning of this article, he specified this mode as  that of {\em seiend-mit}, of being-together-with. For Weyl, remarks such as the one just cited show that  he longed for  a moral quality of being together which we might prefer to call a mode of {\em conviviality} using the later terminology proposed by \cite{Illich:Tools}.\footnote{We can thus express a  distinction to a purely Heideggerian background. That seems   appropriate, because  the person M. Heidegger  was too  closely affiliated with Nazi politics, at least in the first years after their rise to power, to accept that he and his philosophy might have become, without specifications,  a morally attractive point of reference for Weyl in 1948.}
 From Weyl's perspective,  such a mode of convivial being and life might  in turn contribute to the self-understanding of science and become part  of the scientic enterprise, although in this case too he felt only faint ``contours'' of such a hope,  outbalanced by the fear expressed at the end of our last quotation. Nothing assured or assures us of the success of a turn in this direction, however necessary  it  may be. 

\nocite{Weyl:GA}
 \bibliographystyle{apsr}
  \bibliography{a_litfile}

\subsubsection*{ }
\end{document}